\documentclass[11pt]{article}

\RequirePackage[amsthm,amsmath]{imsart}

 \parskip        \baselineskip
\usepackage[mathscr]{eucal}
\usepackage{graphicx}
\usepackage{latexsym}
\usepackage{amssymb}
\usepackage{psfrag}
\usepackage{epsfig}
\DeclareMathAlphabet{\mathpzc}{OT1}{pzc}{m}{it}

\include{psbox}

\usepackage{amsfonts}

\usepackage{graphicx}

\usepackage{amsfonts}

\usepackage{color}

\begin{document}
\bibliographystyle{plain}

\newtheorem{theorem}{\bf Theorem}[section]
\newtheorem{example}{\bf Example}[section]
\newtheorem{definition}{\bf Definition}[section]
\newtheorem{corollary}{\bf Corollary}[section]
\newtheorem{remark}{\bf Remark}[section]
\newtheorem{lemma}{\bf Lemma}[section]
\newtheorem{assumption}{Assumption}
\newtheorem{condition}{\bf Condition}[section]
\newtheorem{proposition}{\bf Proposition}[section]
\newtheorem{definitions}{\bf Definition}[section]
\numberwithin{equation}{section}
\newcommand{\skp}{\vspace{\baselineskip}}
\newcommand{\noi}{\noindent}
\newcommand{\osc}{\mbox{osc}}

\newcommand{\eps}{\varepsilon}
\newcommand{\del}{\delta}
\newcommand{\Rk}{\mathbb{R}^k}
\newcommand{\R}{\mathbb{R}}
\newcommand{\spa}{\vspace{.2in}}
\newcommand{\V}{\mathcal{V}}
\newcommand{\E}{\mathbb{E}}
\newcommand{\I}{\mathbb{I}}
\newcommand{\p}{\mathbb{P}}
\newcommand{\PP}{\mathbb{P}}

\newcommand{\QQ}{\mathbb{Q}}

\newcommand{\lan}{\langle}
\newcommand{\ran}{\rangle}
\def\wh{\widehat}
\newcommand{\defn}{\stackrel{def}{=}}
\newcommand{\txb}{\tau^{\epsilon,x}_{B^c}}
\newcommand{\tyb}{\tau^{\epsilon,y}_{B^c}}
\newcommand{\tilxb}{\tilde{\tau}^\eps_1}
\newcommand{\pxeps}{\mathbb{P}_x^{\eps}}
\newcommand{\non}{\nonumber}
\newcommand{\dist}{\mbox{dist}}

\newcommand{\tilyb}{\tilde{\tau}^\eps_2}
\newcommand{\beq}{\begin{eqnarray*}}
\newcommand{\eeq}{\end{eqnarray*}}
\newcommand{\beqn}{\begin{eqnarray}}
\newcommand{\eeqn}{\end{eqnarray}}
\newcommand{\ink}{\rule{.5\baselineskip}{.55\baselineskip}}

\newcommand{\bt}{\begin{theorem}}
\newcommand{\et}{\end{theorem}}
\newcommand{\deps}{\Delta_{\eps}}

\newcommand{\be}{\begin{equation}}
\newcommand{\ee}{\end{equation}}
\newcommand{\ac}{\mbox{AC}}
\newcommand{\BB}{\mathbb{B}}
\newcommand{\VV}{\mathbb{V}}
\newcommand{\DD}{\mathbb{D}}
\newcommand{\KK}{\mathbb{K}}
\newcommand{\HH}{\mathbb{H}}
\newcommand{\TT}{\mathbb{T}}
\newcommand{\CC}{\mathbb{C}}
\newcommand{\SSS}{\mathbb{S}}
\newcommand{\EE}{\mathbb{E}}

\newcommand{\clg}{\mathcal{G}}
\newcommand{\clb}{\mathcal{B}}
\newcommand{\cls}{\mathcal{S}}
\newcommand{\clc}{\mathcal{C}}
\newcommand{\clv}{\mathcal{V}}
\newcommand{\clu}{\mathcal{U}}
\newcommand{\clr}{\mathcal{R}}
\newcommand{\clt}{\mathcal{T}}
\newcommand{\cll}{\mathcal{L}}

\newcommand{\cli}{\mathcal{I}}
\newcommand{\clp}{\mathcal{P}}
\newcommand{\cla}{\mathcal{A}}
\newcommand{\clf}{\mathcal{F}}
\newcommand{\clh}{\mathcal{H}}
\newcommand{\N}{\mathbb{N}}
\newcommand{\Q}{\mathbb{Q}}

\newcommand{\curvz}{{\bf \mathpzc{z}}}
\newcommand{\curvx}{{\bf \mathpzc{x}}}
\newcommand{\curvi}{{\bf \mathpzc{i}}}
\newcommand{\curvs}{{\bf \mathpzc{s}}}

\newcommand{\tac}{\mbox{\scriptsize{AC}}}
\newcommand{\beginsec}{
\setcounter{lemma}{0} \setcounter{theorem}{0}
\setcounter{corollary}{0} \setcounter{definition}{0}
\setcounter{example}{0} \setcounter{proposition}{0}
\setcounter{condition}{0} \setcounter{assumption}{0}
\setcounter{remark}{0} }

\numberwithin{equation}{section} \numberwithin{lemma}{section}

\begin{frontmatter}
\title{On Uniform Positivity of Transition Densities of Small Noise Constrained Diffusions.}

 \runtitle{Small Noise Asymptotics}

\begin{aug}
\author{Amarjit Budhiraja and Zhen-Qing Chen\\ \ \\
}
\end{aug}

\today

\skp

\begin{abstract}
\noi
Constrained diffusions in convex polyhedral domains with a general oblique reflection field, and with a diffusion coefficient scaled by  a small parameter $\eps > 0$, are considered.
Using an interior Dirichlet heat kernel lower bound estimate for second order elliptic operators in bounded domains from \cite{Zq},
certain uniform in $\eps$ lower bounds on transition densities of such constrained diffusions are established.
These lower bounds together with results from \cite{BiBu2011} give, under additional stability conditions, an exponential leveling property as $\eps \to 0$ for exit times from
suitable bounded domains.

\noi {\bf AMS 2000 subject classifications:} Primary 60F10, 60J60, 60J25; secondary 93E15, 90B15.

\noi {\bf Keywords:} Exponential leveling, reflected diffusions, Dirichlet heat kernel estimates, Skorohod problem, exit time estimates, Friedlin-Wentzell asymptotics.
\end{abstract}

\end{frontmatter}

\section{Introduction}\label{introsec}
\beginsec

Diffusions in polyhedral domains have been extensively studied in the heavy traffic limit theory for stochastic processing networks (see for example, \cite{Rei, harrison, Kus, yamada, BuGhLe}).  In a recent work \cite{BiBu2011} small noise
asymptotics for a general family of such constrained diffusions have been studied. The precise setting there is as follows.
Let $G\subset \Rk $ be convex polyhedral cone with a nonempty interior with the
vertex at origin given as the intersection of half spaces $G_i,
i=1,2,..., N$. Let $n_i$ be the
unit inward normal vector
associated with $G_i$
via the relation
\begin{displaymath}
G_i=\{x\in\Rk \, : \, \langle x,n_i\rangle\geq 0\}.
\end{displaymath}
We will
denote the set $\{x\in\partial G \, : \,\langle x, n_i\rangle=0\}$
by $F_i$.
With each face $F_i$ we associate a unit vector $d_i$ such that
$\lan d_i, n_i\ran >0$. This vector defines the \textit{direction of
constraint} associated with the face $F_i$.  Precise definition of  constrained diffusions considered here is given in Section \ref{prelim}, but
roughly speaking such a process evolves infinitesimally as a diffusion in $\R^k$ and is instantaneously pushed back using the oblique reflection direction
$d_i$ upon reaching the face $F_i$.
Formally, such a process can be represented as a solution of a stochastic integral
equation of the form
\be
\label{ab1009}
X^{\eps}_x(t) = \Gamma\left(x + \int_0^{\cdot} b(X^{\eps}(s)) ds + \eps \int_0^{\cdot} \sigma(X^{\eps}(s)) dW(s)\right)(t),
\ee
where $\Gamma$ is the Skorohod map (see below Definition \ref{defnsp}) taking trajectories with values in $\R^k$ to those with values in $G$, consistent
with the constraint vectors $\{d_i, i = 1, \cdots N\}$.  Under certain regularity assumptions on the Skorohod map (see Condition \ref{cond2.1}) and the
usual Lipschitz conditions on the coefficients $b$ and $\sigma$, the above integral equation has a unique pathwise solution.
One of the main results of \cite{BiBu2011} is an `exponential leveling ' property of exit times from bounded domains for such small noise diffusions.
Such results for diffusions in $\mathbb{R}^k$ have been obtained in \cite{day} which is concerned with asymptotics of Dirichlet problems in bounded domains associated with diffusions with
infinitesimal generator of the form
\begin{equation}\label{genepsdiff}
\cll_{\eps}f(x) =  \frac{\eps^2}{2} \mbox{Tr}(\sigma(x)  D^2f(x)\sigma'(x))+ b(x)\cdot \nabla f(x), \; f \in C_b^2(\R^k) .\end{equation}
The precise result in \cite{day} is as follows. Denote by $Z^{\eps}_x$ the diffusion process governed by the generator $\cll_{\eps}$ and initial distribution $\delta_x$.
Let $B$ be a bounded domain in $\R^k$ and $K$ be an arbitrary compact subset of $B$.  Suppose that all solutions of the ODE $\dot{\xi} = b(\xi)$
with $x =\xi(0)\in B$ converge, without leaving $B$, to a single linearly asymptotically stable critical point. Then, with suitable conditions on the coefficients of the diffusion, for all bounded measurable $f$
$$\sup_{x,y \in K} |\E\left(f(Z^{\eps}_x(\tau^{\eps}_x))\right) - \E\left(f(Z^{\eps}_y(\tau^{\eps}_y))\right)|$$
converges to $0$ at an exponential rate.  Here, $\tau^{\eps}_x = \inf\{t: Z^{\eps}_x(t) \in B^c\}$.  This property is a statement on the long time behavior of the diffusion and
says that although, as $\eps \to 0$, the exit time of the process from the domain approaches $\infty$, the expected values  of functionals of exit location, corresponding to distinct initial conditions, coalesce asymptotically, at an exponential rate. The key ingredient in the proof is the gradient estimate
\be
\label{ab0943}
\sup_{x\in K} |\nabla u^{\eps}(x)| \le c \eps^{-1/2},\ee
where $u^{\eps}$ is the solution of the Dirichlet problem
\beq
 \begin{cases}
     \cll_{\eps}u^{\eps}(x) =0, &  x \in B \\
     u(x) = f(x), &  x \in \partial B .
\end{cases}
\eeq
Diffusions of interest in \cite{BiBu2011} and in the current work are constrained to take values  in domains with
corners and where the constraining mechanism is governed by an oblique reflection field that changes discontinuously from one face of the boundary to another.
To the best of our knowledge there are no regularity (e.g. $C^1$ solutions)  results  known for the associated partial differential equations(PDE)
with oblique reflecting boundary condition.
 In view of this,
 a  probabilistic approach for the study of exponential leveling property for such constrained settings, that `almost' bypasses all PDE estimates
was developed in \cite{BiBu2011}.
The main step in the proof is the construction of certain (uniform in $\eps$) Lyapunov functions under a suitable stability condition which are then used to construct a coupling of the processes
$X^{\eps}_x$, $X^{\eps}_y$ with explicit uniform estimates on exponential moments of time to coupling. The key ingredient
in this coupling construction is a, uniform in $\eps$, minorization condition on transition densities of the
 reflected diffusions (see Condition \ref{regdens}).
The paper \cite{BiBu2011} gave one simple example with constant drift and diffusion coefficients where such a condition is
satisfied. However
the question of when such a minorization
property is available was left as an open problem.

The objective of this work is to answer this question and give general conditions under which Condition \ref{regdens} holds. The main result of the paper is Theorem \ref{mainth}
that shows that under Conditions \ref{cond2.1} and \ref{cond2} the minorization statement in Condition \ref{regdens} holds.
This result, together with Theorem 3.2 in \cite{BiBu2011} then gives general sufficient conditions for an exponential leveling property to hold for a broad family of constrained diffusion processes.   This is noted in Corollary \ref{maincor}.

Once the estimate in Theorem \ref{mainth} is available, the proof of the exponential leveling property does not use any PDE results, however the proof of Theorem \ref{mainth} itself crucially relies on
an interior lower bound estimate for the Dirichlet heat kernel of
 $\mathcal{L}_{\eps}$ over bounded domains
  (\cite{Zq}, see Theorem \ref{T:zhang}).
  In this sense the proofs are not fully probabilistic.

The rest of the paper is organized as follows.  In Section \ref{prelim} we introduce the precise mathematical setting and state our main result (Theorem \ref{mainth}).  In Corollary \ref{maincor} we present
the exponential leveling result that follows on combining Theorem \ref{mainth} with the results in \cite{BiBu2011}.  Finally, in Section \ref{sec:proof} we present the proof of Theorem \ref{mainth}.

The following notation will be used.
Closure, complement, boundary and interior of a subset $B$ of a topological space $S$ will be denoted by $\bar{B}$, $B^c$, $\partial B$ and $B^{\circ}$, respectively.
For a set $B \in \R^k$ and $a \in \R$, we denote by $aB$ the set $\{ax: x \in B\}$.
Borel $\sigma$ field on a metric space $S$ will be denoted as $\clb(S)$.
 Given a metric space $(S,d)$, and subsets $A,B$ of $S$, we will define $\dist(A,B) = \inf_{x\in A, y \in B} d(x,y)$.  If $A = \{x\}$ for some $x \in S$,
we write $\dist(A,B)$ as $\dist(x,B)$. Denote by $C([0,\infty);S)$ (resp. $C([0,T];S)$) the space of continuous functions from $[0, \infty)$ (resp. $[0,T]$)
to a metric space $S$. This space is endowed with the usual local uniform
topology.  For $\eta \in C([0,\infty);\R^k)$, and $t > 0$, define
$|\eta|_t = \sup \sum_{i=1}^n |\eta(t_{i+1}) - \eta(t_i)|,$ where $\sup$ is taken over all partitions $0=t_1< t_2 < \cdots <t_{n+1} =t$.
Lebesgue measure on $\mathbb{R}^k$ will be denoted by $\lambda$.  We denote by $C_b^2(\R^k)$ the space of bounded twice continuously differentiable functions with bounded derivatives.

\section{Main Result}\label{prelim}
\beginsec
 We begin by making precise the constraining mechanism,  in terms of a suitable Skorohod problem, that keeps the diffusion in the polyhedral cone $G$.
Recall the half spaces $G_i$ and vectors $n_i, d_i$, $i = 1, \cdots N$ that were introduced in Section \ref{introsec}..
 For $x\in\partial G$
define
\begin{displaymath}
d(x)=\{v\in\Rk \, : \, v=\sum_{i\in \mbox{In}(x)}\alpha_id_i; \, \alpha_i\geq 0; \, |v|=1\},
\end{displaymath}
where
\begin{displaymath}
\mbox{In}(x)=\{i\in \{1, 2,...,N\}:\,\langle x, n_i\rangle=0\}.
\end{displaymath}
\begin{definition}(Skorokhod Problem).\label{defnsp}
Let $\psi\in C([0,\infty),\R^k)$ be given such that $\psi(0) \in G$.
Then $(\phi, \eta)\in C([0,\infty),G) \times C([0,\infty),\R^k) $
solves the Skorokhod Problem (SP) for $\psi$ (with respect to the data $\{(d_i, n_i), i = 1, \cdots N\}$), if $\phi(0)=\psi(0)$ and if for
all $t\in [0,\infty)$: (1) $\phi(t)=\psi(t)+\eta(t)$, (2) $|\eta|(t)
<\infty$, (3) $|\eta|(t)=\int_0^t\I_{\{\phi(s)\in\partial
G\}}d|\eta|(s)$, and there exists a Borel Measurable function
$\gamma:[0,\infty)\to \R^k$ such that
$\gamma(t)\in d(\phi(t))$, for $d|\eta|$-a.e. $t$ and
$ \eta(t)=\int_0^t\gamma(s)d|\eta|(s)$, $t \ge 0$.
\end{definition}
\noindent Let $C_G([0,\infty);\R^k)$ be the collection of $\psi \in
C([0,\infty);\R^k)$ such that $\psi(0) \in G$.  The domain
$D\subseteq C_G([0,\infty);\R^k)$ on which there is a unique
solution to the Skorohod problem we define the Skorohod map (SM)
$\Gamma$ as $\Gamma(\psi) \doteq \phi$ if $(\phi, \eta)$ is the
unique solution of the Skorohod problem posed by $\psi$. We will
make the following assumption on the regularity of the Skorohod map
defined by the data $\{(d_i,n_i); i=1,2,\ldots, N\}.$
\begin{condition}\label{cond2.1} The Skorohod map is well
defined on all of $C_G([0,\infty);\R^k),$ that is,
$D=C_G([0,\infty);\R^k),$ and the SM is Lipschitz continuous in the
following sense: There exists a constant $K\in(0,\infty)$ such that
for all $\phi_1,\phi_2 \in C_G([0,\infty);\R^k)$:
\begin{equation}\label{sm}
\sup_{0\leq
t<\infty}|\Gamma(\phi_1)(t)-\Gamma(\phi_2)(t)|<K\sup_{0\leq
t<\infty}|\phi_1(t)-\phi_2(t)|.
\end{equation}
\end{condition}
We refer the reader to
\cite{harrison, Dupuis:Ishii:1991, Dupuis:Ramanan:1999}
for sufficient conditions under
which Condition \ref{cond2.1} holds.

Now we introduce the small noise constrained diffusion process that
will be considered in this work. Let $(\Omega, \mathcal{F}, \p)$ be
a complete probability space on which is given a filtration
$\{\mathcal{F}_t\}_{t\geq 0}$ satisfying the usual hypothesis. Let
$(W(t), \mathcal{F}_t)$ be a $k$-dimensional standard Wiener process
on the above probability space.  Let
 $\sigma : G\to \mathbb{R}^{k\times k}$ , $b: G\to \Rk$
be mappings satisfying the following condition.
\begin{condition}\label{cond2}
(i) There exists $\gamma_1\in (0, \infty)$ such that
 \be  \label{cond1}
  |b(x)- b(y)|\leq \gamma_1 |x-y| \ \hbox{ and } \   |b(x)| \le \gamma_1
  \quad \hbox{for every } x, y\in G
\ee
(ii) There exists $\gamma_2\in (0, \infty)$ such that
 \be
 |\sigma(x)-\sigma(y)| \leq \gamma_2 |x-y|
\  \hbox{ and } \ |\sigma(x)|  \leq \gamma_2
  \quad \hbox{for } x, y\in G.   \label{cond2eq} \ee
(iii) There exists $\underline{\sigma} \in (0, \infty)$ such that for all $x \in G$ and $v \in \mathbb{R}^{k}$
$$v' \sigma(x) \sigma(x)' v \ge \underline{\sigma}
|v|^2. $$
\end{condition}
Given $\eps > 0$, let $X^{\eps}$ be
the unique strong solution of the following stochastic integral
equation: \be  X(t)=\Gamma\left(x + \int_0^. b(X(s))ds +
\eps\int_0^. \sigma(X(s))dW(s)\right)(t),\; t \ge 0,\label{eqn1} \ee
Existence of strong solutions and pathwise uniqueness for \eqref{eqn1} is a consequence of the Lipschitz property
of the coefficients and of the Skorohod map (see \cite{Dupuis:Ishii:1991}).  It is convenient to have the process $X$ for various initial
conditions and values of $\epsilon$ to be defined on a common canonical space.  Indeed, one can find a filtered measurable space, which we denote again as  $(\Omega, \mathcal{F},  \{\mathcal{F}_t\})$, on which is given a family of
 probability
measures, $\{\mathbb{P}^{\eps}_x \}_{x\in G}$,  and continuous adapted stochastic processes $Z$, $Y$ and $W$ such that for all $x \in G$, under
$\mathbb{P}^{\eps}_x$,  $\{W(t ), \{\mathcal{F}_t \}_{t \ge0} \}$ is a $k$ -dimensional standard Wiener process and $(Z, W,Y )$ satisfy
$\mathbb{P}^{\eps}_x$ a.s. the integral equation
\beqn Z(t)&=& \Gamma\left(x + \int_0^. b(Z(s))ds +
\eps\int_0^. \sigma(Z(s))dW(s)\right)(t),\nonumber\\
&=& x + \int_0^t b(Z(s))ds +
\eps\int_0^t \sigma(Z(s))dW(s) + \mathbb{D}Y(t), \; t \ge 0,\label{ab1719}
\eeqn
where $\mathbb{D}= (d_1, \cdots d_N)$.
 Also, for every $\eps > 0$,  $(Z, \{\PP^{\eps}_x\}_{x\in G})$ is a strong Markov family (cf. \cite{amarjitlee}).

It can be shown that for $t > 0$ and $x \in G$, the measure $\PP_x^{\eps}(Z(t) \in \cdot)$ is absolutely continuous with respect to the Lebesgue measure on $G$ (see Lemma 5.7
in \cite{amarjitlee}). Denote by $p_{\eps}(t, x,y)$ the probability density of $Z(t)$ under $\PP^{\eps}_x$, namely
$$\PP^{\eps}_x(Z(t) \in A) = \int_A p_{\eps}(t, x,y) dy ,\; A \in \mathcal{B}(G).$$
Let, for $r > 0$ and $x_0 \in G$, $\BB_r(x_0) = \{x \in G: |x-x_0| \le r\}$.  When $x_0=0$, we simply write $\BB_r$.

The following minorization condition was introduced in \cite{BiBu2011} (see Condition 3.1 therein).
\begin{condition}
	\label{regdens}
	For every $t_1>0$, $\eps_0 > 0$ and $0 < M_1 < M < \infty$, there exists a Borel $E \subset \BB_{M_1}$ with $\lambda(E) > 0$ and  $\kappa > 0$ such that, for all $\eps \in (0, \eps_0)$,
	$$\eps^{2k} p_{\eps}(t_1\eps^2, x,z) \ge \kappa
 \quad \hbox{for all } x \in \eps^2\BB_M \mbox{ and } \lambda \mbox{-a.e. } z \in \eps^2E.$$
\end{condition}
The above condition played a key role in \cite{BiBu2011} in proving an exponential leveling result for exit times from bounded domains which we now describe.
Let
\begin{displaymath}
\mathcal{C}=\left\{-\sum_{i=1}^N\alpha_id_i : \, \alpha_i\geq 0;\, i\in
\{1,...,N\}\right\} .
\end{displaymath}
The paper \cite{ramiamarjit} shows that under Conditions \ref{cond2.1} and \ref{cond2} and Condition \ref{cond3} below  the constrained diffusion $X^{\eps}$ is positive recurrent.
\begin{condition}
	\label{cond3}
For some $\del > 0$
	\be
	\label{ab1120}
	b(x) \in \clc(\delta) = \{v\in \mathcal{C}
	\, : \, \mbox{dist}(v, \partial\mathcal{C})\geq\delta\}
\quad \hbox{for all } x \in G.
	\ee
\end{condition}
Let $B$ be a bounded open subset of $G$.  Suppose that
$0 \in B$ and $\partial B = \partial \bar{B}$.

For $x \in G$ we denote by $\xi_x$ the unique solution of the integral equation
\be \label{ab1227}
\xi_x(t) = \Gamma \left( x + \int_0^{\cdot}b(\xi_x(s)) ds \right)(t),
\quad t \ge 0 .
\ee
Also, let $\cls_x = \{\xi_x(t), t \ge 0\}$.  For $\gamma > 0$, let
$B_{\gamma} = \{x \in B: \dist(\cls_x, \partial B) \ge \gamma\}$ and let $B_0 = \bigcup_{\gamma > 0} B_{\gamma}$.
  Clearly $B_0 \subset B$ and under Condition \ref{cond3} one can show that for every $\gamma$ sufficiently small, there is a $r > 0$ such that
$\BB_{r} \subset B_{\gamma}$ and
for every $x \in B_0$, $\xi_x(t)$ converges to $0$, as $t\to \infty$, without
leaving $B$ (cf. Lemma 2.1 in \cite{BiBu2011}).

Denote by $\clh$ the collection of all $\psi:\R^+ \to \R^+$ satisfying $\sup_{t >0} \frac{|\psi(t)|}{(1+ \log^+(t))^{q}} < \infty$, for some
$0 < q < 1$,
and
 $\sup_{r>0} \frac{\sup_{s,t \ge 0: |t-s| \le r}|\psi(t)-\psi(s)|}{r^m +1} < \infty $
 for some $m \ge 1.$  Let $\tau \equiv \tau(B) = \inf \{t \ge 0: Z(t) \in B^c\}$.

The following is the exponential leveling result from \cite{BiBu2011}.
\bt \label{levelmain}
Suppose that Conditions \ref{cond2.1}, \ref{cond2},
\ref{regdens} and \ref{cond3} holds.
Let $\KK$ be a compact subset of $B_0$. Then the following hold.\\
\noindent (i)
 For every bounded measurable map $f: \bar B \to \R$ there is a $\del_1 \in (0, \infty)$ such that
$$
\lim_{\eps \to 0}\sup_{x,y \in \KK}e^{\del_1/\eps} |\E_x^{\eps} f(Z(\tau)) - \E_y^{\eps}f(Z(\tau))| = 0.
$$
\noindent (ii) For every  $\psi\in \clh$
$$\limsup_{\eps \to 0} \sup_{x,y \in \KK}|\E_x^{\eps} \psi(\tau) - \E_y^{\eps}\psi(\tau)| < \infty.$$
\et
The main result of the current work says that Condition \ref{regdens} in Theorem \ref{levelmain} can be dropped.  Specifically, we prove the following result.

\medskip

\begin{theorem}
	\label{mainth}
	Suppose Conditions \ref{cond2.1} and \ref{cond2} hold. Then Condition \ref{regdens} holds.
	\end{theorem}
	As an immediate consequence of this result we have the following.
	\begin{corollary}
		\label{maincor}
		Suppose that Conditions \ref{cond2.1}, \ref{cond2} and  \ref{cond3} hold.  Then properties (i) and (ii) in the statement of Theorem \ref{levelmain}
		hold.
\end{corollary}

Rest of this paper is devoted to the proof of Theorem \ref{mainth}.

\section{Proof of Theorem \ref{mainth}.}
\label{sec:proof}
For $x \in \R^k$, denote by $Q_{\cll,x}$ the probability law of the diffusion with initial distribution $\delta_x$ and infinitesimal generator $\cll$:
\begin{equation}\cll f(x) = b(x)\cdot \nabla f(x) + \frac{1}{2} \mbox{Tr}(\sigma(x) D^2f(x)\sigma'(x)), \; f \in C_b^2(\R^k),\label{eq:eq115}
\end{equation}
where $b$, $\sigma$ satisfy Condition \ref{cond2}
and $\sigma' (x)$ is the transpose of the matrix $\sigma(x)$.
Denote by $\xi$ the canonical coordinate process on $C([0,\infty);\R^k)$ and given
 any bounded $C^2$-smooth open set $B$ in $\R^k$, let
$$\tau_B = \inf\{t\ge 0: \xi(t) \in B^c\}.$$
Then,
for every $t > 0$ and $\eps > 0$, there is a jointly continuous
function
 $p^{\dag}_{\cll,B}(t, \cdot, \cdot): B \times B \to \R_+$
such that for all $x \in B$ and Borel $A \subset B$
$$Q_{\cll,x}(\xi(t) \in A; t<\tau_B) = \int_A p^{\dag}_{\cll,B}(t,x,y) dy.$$
The function $p^{\dag}_{\cll,B}(t,x,y)$ is called Dirichlet heat kernel
for $\cll$ in $B$.

Also, define $\wh \cll_{\eps}$ by replacing $\sigma$
and $b$ in \eqref{eq:eq115} with $b_{\eps}(\cdot) = b(\eps^2\cdot)$ and $\sigma$ with $\sigma_{\eps}(\cdot) = \sigma(\eps^2\cdot)$.
Let $a(x)=(a_{ij}(x))=\sigma (x) \sigma' (x)$ and $\wh b (x)=(\wh b_1(x),
\dots, \wh b_k(x))$ with $\wh b_k(x)=\sum_{j=1}^k \frac{\partial a_{kj}(x)}
{\partial x_j}$.
Note that $\wh b$ is defined as a bounded function almost everywhere
on $\R^k$, as each $a_{ij}$ is a bounded Lipschitz function.
We can rewrite $\wh \cll_{\eps}$ almost everywhere on $\R^k$ as follows:
$$
\wh \cll_{\eps}=\frac12 \nabla (a(\eps^2 x)\nabla )
 + (b(\eps^2 x) -\eps^2 \wh b(\eps^2 x))\nabla.
$$
Various heat kernel estimates for operators of the above type
have been studied in literature.
 
The following interior lower bound estimate follows immediately
 from \cite[Lemma 3.3]{Zq} and a finite covering argument.
(A sharp two-sided Dirichlet heat kernel estimate is available
from \cite[Theorem 2.1]{Ria}, but we do not need it in this paper.
Though it is assumed in \cite{Zq} that the dimension
$k\geq 3$, the interior lower
bound estimate \cite[Lemma 3.3]{Zq} in fact holds for any $k\geq 1$;
see \cite{Ria}.)

\begin{theorem}
	\label{T:zhang}
Let $B=B(x_0, R)$ be an open ball in $\R^k$.
For each given $\eps_0>0$, $t>0$ and $0<\gamma <1$, there is a constant $c>0$ so that
$p^{\dag}_{\wh\cll_{\eps},B}(t,x,y) \ge c$ for every $x, y\in B(x_0, \gamma R)$ and $\eps \in (0, \eps_0]$.	
\end{theorem}

\subsection{Proof of Theorem \ref{mainth}.}
Fix $t_1,M , \eps_0 > 0$ and $M_1 \in (0,M)$.  Also, fix $x_0 \in G^{\circ}$ such that $|x_0| < M_1$.  Choose $r_2 > 0$ such that
$\BB_{r_2}(x_0)\subset G^{\circ} \cap \{x: |x| < M_1\}$.  Fix $0 < r_0 < r_1 < r_2$.  Also, fix $t_2, t_3 > 0$ such that $t_2+t_3 = t_1$.
We will now show that Condition \ref{regdens} holds  with $E = \BB_{r_0}(x_0)$.
Let $\varphi: G \to [0,1]$ be a continuous function such that $\varphi(x) = 1$ for all $x \in \BB_{r_0}(x_0)$ and
$\varphi(x) = 0$ for all $x \in (\BB_{r_1}(x_0))^c$.  We first show that
\begin{equation}
	\label{eq:eq749}
	\mbox{ For some } \kappa_0 > 0, \
\int_{\R^k} \varphi(z/\eps^2) p_{\eps}(t_2\eps^2, x, z) dz \ge \kappa_0
 \quad \hbox{for all } x \in \BB_{M\eps^2} \mbox{ and }
 \eps \in (0, \eps_0]
\end{equation}
For this, note that
$$\int_{\R^k} \varphi(z/\eps^2) p_{\eps}(t_2\eps^2, x, z) dz = \EE_x^{\eps}\left[\varphi\left(\frac{Z(t_2\eps^2)}{\eps^2}\right)\right]$$
and,  $\PP_x^{\eps}$ a.s.,
$$Z(t) = \Gamma\left(x + \int_0^. b(Z(s))ds +
\eps\int_0^. \sigma(Z(s))dW(s)\right)(t).$$
Therefore, letting $Z^{\eps}(t) = \frac{Z(t_2\eps^2)}{\eps^2}$ and using standard scaling properties of the Skorohod map
$$Z^{\eps}(t)  = \Gamma\left(\frac{x}{\eps^2} + \int_0^. b(\eps^2Z^{\eps}(s))ds +
\int_0^. \sigma(\eps^2 Z^{\eps}(s))dW^{\eps}(s)\right)(t),$$
$\PP_x^{\eps}$ a.s., where $W^{\eps}(t) = \frac{1}{\eps}W(t\eps^2)$.
Suppose now that \eqref{eq:eq749} fails.  Then there exists a sequence $\eps_n \to 0$ and $\{x_n\} \subset G$ such that $x_n \in \BB_{M\eps_n^2}$ and
$$\lim_{n \to \infty} \EE_{x_n}^{\eps_n}\left[ \varphi(Z^{\eps_n}(t_2))\right] = 0.$$
Without loss of generality we can assume that $x_n/\eps_n^2$ converges to some $\bar{x} \in \BB_M$. A standard weak convergence argument now shows that
\begin{equation}
	\label{eq:eq802}
0 = \lim_{n\to \infty} \EE_{x_n}^{\eps_n}\left[ \varphi(Z^{\eps_n}(t_2))\right] = \tilde E\left [\varphi(\tilde X(t_2))\right],
\end{equation}
where
$\tilde X$ is a $G$ valued continuous stochastic process given on some  probability space $(\tilde \Omega, \tilde \clf, \tilde P)$,
supporting a $k$ dimensional Brownian motion $\tilde W$, such that
$$
\tilde X(t) = \Gamma \left (\bar{x} + b(0)\iota + \sigma(0)\tilde W\right)(t), \; t \ge 0.$$
From Condition \ref{cond2} (iii) it follows that (cf. Lemma 5.7\cite{amarjitlee})
$$\tilde \E\left [\varphi(\tilde X(t_2))\right] \ge \tilde \PP(\tilde X(t_2) \in \BB_{r_0}(\bar{x})) > 0.$$
This contradicts \eqref{eq:eq802} and therefore \eqref{eq:eq749} follows.

Note that if $Z$ is a diffusion process on $\R^k$ having transition density
function $p(t, x, y)$ with respect to the Lebesgue measure on $\R^k$,
then $Y(\cdot)=\eps^{-2} Z(\eps^2 \cdot)$ has transition density function
$p^Y(t, x, y)=\eps^{2k}p(\eps^2t, \eps^2 x, \eps^2 y)$.

Let $\cll_{\eps}$ be as introduced in \eqref{genepsdiff}.
Its associated diffusion is given by
$$ dZ(t)= \eps \sigma (Z(t))dW(t) + b(Z(t)) dt.
$$
Using Brownian scaling, we see that $Y(t)=\eps^{-2} Z(\eps^2 t)$ satisfies
$$ dY(t) = \sigma ( \eps^{2} Y(t))dW(t) + b(\eps^{2} Y(t)) dt,
$$
and thus has infinitesimal generator $\wh \cll_{\eps}$ introduced above Theorem \ref{T:zhang}, namely,
$$ \wh \cll_{\eps} f(x) =  \frac{1}{2} \mbox{Tr}
\left( \sigma(\eps^{2} x)  D^2 f(x) \sigma'(\eps^{2} x) \right)
+ b(\eps^{2} x)\cdot \nabla f(x), \; f \in C_b^2(\R^k).
$$

Next, let $p^{\dag}_{\eps} = p^{\dag}_{\cll_{\eps}, B_{\eps}}$, where $B_{\eps} = \eps^2(\BB_{r_2}(x_0))^{\circ}$. Then by the above scaling relation,
we conclude that
$$
\eps^{2k}p^{\dag}_{\eps}(\eps^2t, \eps^2 x, \eps^2 y) = p^{\dag}_{\wh\cll_{\eps}, (\BB_{r_2}(x_0))^{\circ}}(t, x, y)
$$
for every $t>0$ and $x, y \in (\BB_{r_2}(x_0))^{\circ}$.

It follows from Theorem \ref{T:zhang}
that there is a constant $\kappa_1>0$ so that for every
$\eps \in (0, \eps_0]$
and for $y \in \eps^2
\BB_{r_1}(x_0)$ and $z \in \eps^2 \BB_{r_1}(x_0)$,
\begin{align}
	\eps^{2k} p^{\dag}_{\eps}(\eps^2t_3, y, z) =& p^{\dag}_{\wh\cll_{\eps}, (\BB_{r_2}(x_0))^{\circ}}(t_3,  \eps^{-2}y, \eps^{-2}z)
\geq \kappa_1 . \label{eq:eq820}
\end{align}
Finally, take $E= \BB_{r_1}(x_0)$.
Then
for every  for $x \in \BB_{M\eps^2}$ and  $z \in \eps^2E$,
\begin{align*}
	\eps^{2k}p_{\eps}(t_1\eps^2, x,z) \ge & \eps^{2k} \int_{\eps^2\BB_{r_1}(x_0)} p_{\eps}(t_2\eps^2, x,y)p_{\eps}(t_3\eps^2, y,z) dy\\
	\ge & \left(\inf_{y, z \in \eps^2\BB_{r_1}(x_0)}
\eps^{2k}p_{\eps}^{\dag}(t_3\eps^2, y,z)\right) \int_{\eps^2\BB_{r_1}(x_0)} p_{\eps}(t_2\eps^2, x,y) dy\\
	\ge & \kappa_1 \int_{\R^k}\varphi(\frac{y}{\eps^2})p_{\eps}(t_2\eps^2, x,y) dy\\
	\ge & \kappa_1 \cdot \kappa_0,
\end{align*}
where the third inequality uses \eqref{eq:eq820} and the fourth follows from \eqref{eq:eq749}. The result now follows on taking $\kappa = \kappa_1 \cdot \kappa_0$.
\qed

 \noindent{\bf Acknowledgement:} Research
   of the first author is supported in
 part by the National Science Foundation (DMS-1004418,DMS-1016441, DMS-1305120) and the Army Research
 Office (W911NF-10-1-0158).
 Research
   of the second author is supported in
 part by  NSF Grant   DMS-1206276  and NNSFC Grant 11128101.


\skp

{\sc

\bigskip\noi
A. Budhiraja\\
Department of Statistics and Operations Research\\
University of North Carolina\\
Chapel Hill, NC 27599, USA\\
email: budhiraj@email.unc.edu

\skp

\noi
Z.-Q. Chen\\
Department of Mathematics\\
University of Washington\\
Seattle, WA 98195, USA\\
email: zqchen@uw.edu

}

\end{document}